\documentclass{article}
\usepackage{amsmath, lscape, graphicx, enumerate, hyperref}
\title{\textbf{On M\=adhava and his correction terms for the M\=adhava-Leibniz series for $\pi$}}
\author{\textbf{V. N. Krishnachandran}\\[1mm]
Sreepathi Institute of Management and technology\\
Vavanoor - 679533, Palakkad, Kerala, India \\
(Email: {\tt krishnachandranvn@gmail.com})}
\date{}
\allowdisplaybreaks
\begin{document}
\maketitle
\tableofcontents
\newpage
\begin{abstract}
This paper is intended to serve two purposes: one, to present an account of the life of Sangamagr\=ama M\=adhava, the founder of the Kerala school of astronomy and mathematics which flourished during the 15th - 18th centuries, based on modern historical scholarship and two, to present a critical study of the three enigmatic correction terms, attributed to M\=adhava, for obtaining more accurate values of $\pi$ while computing its value using the M\=adhava-Leibniz series. For the second purpose, we have collected together the original Sanskrit verses describing the correction terms, their English translations and their presentations in modern notations. The Kerala rationale for these correction terms are also critically examined. The general conclusion in this regard is that, even though the correction terms give high precision approximations to the value of $\pi$,  the rationale presented by Kerala authors is not strong enough to convince modern mathematical scholarship.  

The author has extended M\=adhava's results by presenting  higher order correction terms which yield better approximations to $\pi$ than the correction terms attributed to M\=adhava. The various infinite series representations of $\pi$ obtained by M\=adhava and his disciples from the basic M\=adhava-Leibniz series  using M\=adhava's  correction terms are also discussed. A few more such series representations using the better correction terms developed by the author are also presented.  The various conjectures regarding how M\=adhava might have originally arrived at the correction terms are also discussed in the paper. 
\end{abstract}
\qquad\\[1mm]
\section{Introduction}
The achievements of the mathematician-astronomers of the Kerala school of astronomy and mathematics which  flourished during the 14th-17th centuries CE were first brought to the attention of the western mathematical scholarship with the publication of a now famous paper title ``On the Hindu Quadrature of the Circle, and the Infinite Series of the Proportion of the Circumference to the Diameter Exhibited in the Four S\=astras, the Tantra Sangraham, Yucti Bh\=ash\=a, Carana Padhati, and Sadratnam\=ala'' by C. M. Whish in the Transactions of the Royal Asiatic Society of Great Britain and Ireland in 1834 (see \cite{Whish}). (For a fascinationg account of the life and work  of Whish see \cite{MacTutor3}.) These achievements are now well understood and are widely known thanks to the efforts of scholars more than a century after the publication of Whish's paper(see, for example, \cite{Jush}, \cite{MacTutor2}, \cite{Plofker} Ch. 7, \cite{Uta} p.202, \cite{Victor}).  

Among the mathematicians of the Kerala school, the name of Sangamagr\=ama M\=adhava, who is the considered the founder of this school, stands out prominently even though direct knowledge about his life and work is very patchy and scanty. M\=adhava is well-known as the inventor the series expansions of the trigonometric sine and cosine functions and he is also well-known for providing proofs of these series expansions which meet current standards of mathematical rigour. Among these series, as a special case of the series for the arc-tangent function, M\=adhava has given an infinite series for evaluating the mathematical constant $\pi$. This series was later rediscovered in Europe by Leibniz in around 1675 (see \cite{Roy2}). What is amazing is that M\=adhava has also given three different successively better approximations to the error in computing the value of $\pi$ using the series by truncating it to get a finite computable expression. These approximations generally referred to as  ``correction terms'' yield remarkably accurate values to $\pi$. But what is really unfortunate is that the extant works the Kerala school of mathematics give only rationales for two of these correction terms but do not provide a rationale for the third and the best approximation. This prompted historians of mathematics to speculate on the possible method by which these approximations might have been derived. This paper presents a critical summary of the various approaches that have been adopted to unravel the mystery of the correction terms and also how these correction terms have been put to use to derive infinite series representations of $\pi$ which converge faster than the M\=adhava-Leibniz series. 

The rest of the paper is organised as follows.  Section 2 of the paper presents a new perspective on the life and times of M\=adhava which is at variance with the generally accepted view on this matter. This perspective is based on an elaborately argued exposition of the issues in a book on history of mathematics in India published in 2008 (see \cite{Divakaran}). In the next two sections we present  the M\=adhava-Leibniz series and the correction terms both in modern notations and in the form they have been presented by mathematicians of the Kerala school.  Section 5 is devoted to introducing the notion of a measure of inaccuracy of  a correction term. Section 6 discusses at length the Kerala rationale for the first two correction terms and also, in the spirit of the Kerala rationale, our rationale  for the third correction term. Some properties of the measure of inaccuracy are considered in Section 7. In Section 8, more infinite series representations of $\pi$ are considered all of which appear in the mathematical texts of the Kerala school and all of which converge much faster than the basic M\=adhava-Leibniz series. In the next section we consider the various conjectures regarding how M\=adhava might have obtained the correction terms. The paper concludes by raising the question whether M\=adhava actually used these series and the various correction terms to compute the numerical value of $\pi$.

\section{Sangamagr\=ama M\=adhava: A new perspective on his life}
%
%
Practically very little is known with absolute certainty about Sangamagr\=ama M\=adhava the person, his life and his times. 
The commonly accepted view among scholars has always been that Sangamagr\=ama is the town of Irinj\=alakku\d da, about 60 km north of Cochin and some 70 km south of the Ni\d la river (see \cite{Divakaran} p.286, \cite{History} p.51). Other than some local legends (see, for example, \cite{Ulloor} p.104) and some highly stretched interpretations of place names and the connection between the house name of M\=adhava and the house name of a local Namp\=utiri family there appears to be no other concrete evidence to connect M\=adhava with Irinj\=alakku\d da. 

From scattered references to M\=adhava found in diverse manuscripts, historians of Kerala mathematics have pieced together some bits of information. In a manuscript preserved in the Oriental Institute, Baroda, M\=adhava has been referred to as {\em M\=dhavan v\=e\d{n}v\=r\=h\=d\=n\-m kart\=a ... M\=adhavan Ilaññippa\d{ll} Empr\=an} (see \cite{History} p.51). Firstly, it has been noted that the epithet ``{\em Empr\=an}'' in the name is a reference to a certain community he would have belonged to. ``{\em Empr\=an}'' is a shortened form of the word ``{\em Empr\=antiri}" and an {\em Empr\=antiri} is a member of the community of Brahmins, considered somewhat inferior in status to the Namp\=utiris (see \cite{Plofker} p.218), who have migrated from Tulu Nadu to Kerala and settled in Kerala. Tulu Nadu is not an administrative unit and is   generally identified as the region in the southwestern coast of India consisting of the Dakshina Kannada and Udupi districts of Karnataka state and Kasaragod district of Kerala state. Most members of the community speak Tulu language.

The term ``{\em Ilaññippa\d l\d li}" has been identified as a reference to the house name of M\=adhava. This is corroborated by M\=adhava himself. In his short work on the moon’s positions titled {\em Ve\d nv\=ar\=oha}, M\=dhava says that he was born in a house named {\em baku\d l\=adhi\d s\d tita . . . vih\=ara} (see \cite{Sphuta} p.12). This is clearly a Sanskritisation of {\em Ilaññippa\d{ll}i}. {\em Ilaññi} is the Malayalam name of the evergreen tree Mimusops elengi and the Sanskrit name for the same is {\em Bakuļa}. The word {\em pa\d{ll}i} had been in use to refer to a Buddhist retreat house and the word continued to be used as suffixes to names of houses and places in Kerala even after the disappearance of Budhism from Kerala. The Sanskrit equivalent of {\em pa\d{ll}i} is {\em vihāra}. {\em baku\d lādhi\d s\d thita} can be translated as ``occupied or inhabited by {\em bakula}". 

Incidentally, the Sanskrit house name {\em baku\d lādhi\d s\d thita ... vihāra} has also been interpreted as a reference to the Malayalam house name {\em Iraññi ninna pa\d{ll}i} and some historians have tried to identify it with one of two currently existing houses with names {\em Iriññanava\d{ll}i} and {\em Iriññātapa\d{ll}i} both of which are located near Irinjalakuda town in central Kerala. This identification is far fetched because both names have neither phonetic similarity nor semantic equivalence to the word ``{\em Ilaññippa\d{ll}i}".

Most of the writers of astronomical and mathematical works who lived after M\d adhava's period have referred to M\d adhava as `'{\em Sangamagr\=ama M\=adhava}" and as such it is important that the real import of the word ``Sangamagr\=ama" be made clear. The general view among many scholars is that Sangamagr\=ama is the town of Irinjalakkuda some 70 km south of the Nila river and about 70 km south of Cochin. It seems that there is not much concrete ground for this belief except perhaps the fact that the presiding deity of an early medieval temple in the town, the Koodalmanikyam Temple, is worshipped as Sangameswara meaning the Lord of the Samgama and so Samgamagrama can be interpreted as the village of Samgameswara. But there are several places in Karnataka with {\em samgama} or its Dravidian equivalent {\em kū\d dala} in their names and with a temple dedicated to Samgamesvara, the lord of the confluence. (Kudalasangama in Bagalkote district in Karnataka State, which lies at the point of confluence of the rivers Krishna and Malaprabha, is one such place with a celebrated temple dedicated to the Lord of the Samgama.)

Interestingly, there is a small town on the southern banks of the Nila river, around 10 km upstream from Tirunavaya, called K\=u\d dall\=ur. The exact literal Sanskrit translation of this place name is Sangamagrama: {\em kū\d tal} in Malayalam means a confluence (which in Sanskrit is {\em samgama}) and {\em ūr} means a village (which in Sanskrit is {\em gr\d ama}). Also the place is at the confluence of the Nila river and its most important tributary, namely, the Kunti river. (There is no confluence of rivers near Irinjalakuada.) Incidentally there is a still existing Namp\= utiri (Malayali Brahmin) family by name Kūtallūr Mana a few kms away from the Kudallur village. The family has its origins in Kudallur village itself. For many generations this family hosted a great Gurukulam specialising in Vedanga. Perhaps it is not coincidence that  the only available manuscript of {\em Sphu\d tacandr\=apti}, a book authored by Madhava, was obtained from the manuscript collection of Kūtallūr Mana. This may not be purely accidental; indeed, M\=adhava might have had some association with Kūtallūr Mana.

Thus the most plausible possibility is that the forefathers of M\=adhava migrated from the Tulu land or thereabouts to settle in K\=udall\=ur village, which is situated on the southern banks of the Nila river not far from Tirunnavaya, a generation or two before his birth and lived in a house known as {\em Ilaññippa\d{ll}i} whose present identity is unknown.

Regarding M\=adhava's times there are no definite evidences to pinpoint the period during which M\=adhava lived. In his {\em Venv\=aroha}, M\=adhava gives a date in 1400 CE as the epoch. M\=adhava's pupil Paramesvara, the only known direct pupil of M\=adhava, is known to have completed his seminal work {\em Drigganita} in 1430 and Paramesvara's date has been determined as c.1360-1455. From such circumstantial evidences historians of Kerala mathematics have assigned the date c.1340 - c.1425 to M\=adhava. 

\section{The M\=adhava-Leibniz series and the correction terms}
The following is the M\=adhava-Leibniz series for $\pi$:
\begin{equation}\label{eq1a}
\frac{\pi}{4} = 1 - \frac{1}{3}+\frac{1}{5} - \frac{1}{7} + \cdots + (-1)^{n-1}\frac{1}{2n-1}+\cdots
\end{equation}
Truncating the series after the $n$-th term we have the following approximate expression for $\pi$:
\begin{equation}\label{eq2}
\frac{\pi}{4} \approx 1 - \frac{1}{3}+\frac{1}{5} - \frac{1}{7} + \cdots + (-1)^{n-1}\frac{1}{2n-1}.
\end{equation}
Denoting by $F(n)$ the absolute value of the error in the approximate expression for $\pi$, we may write
\begin{equation}\label{eq1}
\frac{\pi}{4} = 1 - \frac{1}{3}+\frac{1}{5} - \frac{1}{7} + \cdots + (-1)^{n-1}\frac{1}{2n-1}+ (-1)^nF(n).
\end{equation}
The function $F(n)$ as defined in Eq.\eqref{eq1} is called the ``correction term'' that should be added to or subtracted from the right hand side of Eq.\eqref{eq2} to get the exact value of $\frac{\pi}{4}$. Note that $F(n)$ is the correction term after the $n$-th term in the series in \eqref{eq1}.

The Kerala mathematical literature gives the following three approximations for the correction term $F(n)$, all attributed to M\=adhava, each giving successively better approximations than the previous one:
\begin{align}
F_1(n) & = \frac{1}{4n}\\
F_2(n) & = \frac{n}{4n^2+1}\\
F_3(n) & = \frac{n^2+1}{n(4n^2+5)}
\end{align}
\section{The series and the correction terms in Kerala mathematical literature}\label{S4}
The surviving texts containing results on infinite series for $\pi$ are {\em Tantrasangraha} by Ke\d lallur N\=ilakan\d{t}ha Somay\=aj\=i (1444 - 1545), a commentary on it called {\em Yuktid\=ipika-Laghuvivrtti} by \'Sankara V\=ariyar (c.1500 - c.1560), a work in Malayalam titled  {\em Yuktibha\d sa} by
Jye\d sthadeva (c.1500 - c.1575) and the work {\em Kriy\=akramakari}  by \'Sankara V\=ariyar and his student
Mahi\d samangalam N\=arayana. All these works are in Sanskrit except the {\em Yuktibh\=a\d sa} whic is written in Malayalam, the language of Kerala. These works provide a summary of major
results on infinite series discovered by Kerala mathematicians.

The Sanskrit verses describing the correction terms and their English translations given below are as quoted in \cite{Raju} pp.173-174. For another wording of the translations one may consult \cite{Rajagopal}. 
\subsection{First correction term}
The Kerala texts do not explicitly state the first correction term $F_1(n)$ (see Section \ref{Rationale} and Eq.\eqref{eq17}). It is used as a preliminary step in the formulation of the more accurate correction term $F_2(n)$. 

\subsection{Second correction term}
The {\em Yuktid\=ipika-Laghuvivrtti} (Chapter 2: Verses 271 - 274) commentary of Tantrasangraha presents the second correction term in the following verses:
\begin{center}
\includegraphics[width=8cm]{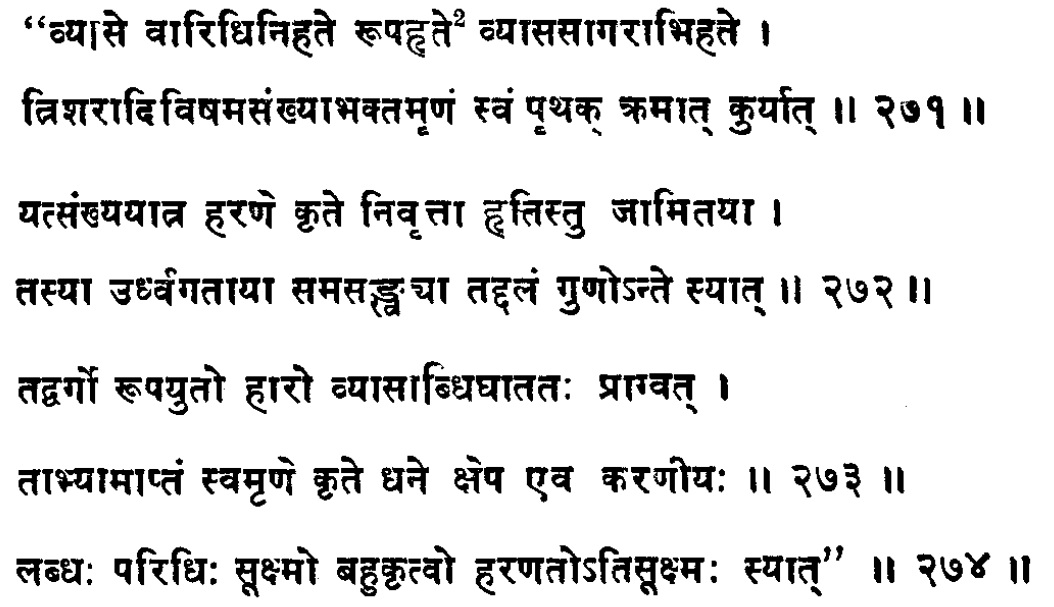}
\end{center}
Here is an English translation of the verses: 
\begin{quote}
``To the diameter multiplied by 4 alternately add and subtract in order the diam­eter multiplied by 4 and divided separately by the odd numbers 3, 5, etc. That odd number at which this process ends, four times the diameter should be multiplied by the next even number, halved and [then] divided by one added to that [even] number squared. The result is to be added or subtracted according as the last term was subtracted or added. This gives the circumference more accurately than would be obtained by going on with that process.''
\end{quote}
Taking the diameter of a given circle as $d$, the result stated in the verses can be presented as follows:
\begin{equation}
\text{Circumference} = 4d - \frac{4d}{3} + \frac{4d}{5} - \cdots \pm \frac{4d}{p}\mp 4d\times  \frac{\left(\frac{p+1}{2}\right)}{(p+1)^2 + 1} 
\end{equation}
The correction term is 
\begin{align*}
f_2(p)&=\frac{\left(\frac{p+1}{2}\right)}{(p+1)^2 + 1}\\
&=\frac{p+1}{2(p^2+2p+2)}.
\end{align*}
If we let $p=2n-1$ this yields $F_2(n)$.
\subsection{Third correction term}
The {\em Yuktid\=ipika-Laghuvivrtti} commentary of Tantrasangraha also gives the correction term $f_3(n)$  in the following verses (Chapter 2: Verses 295 - 296): 
\begin{center}
\includegraphics[width=8cm]{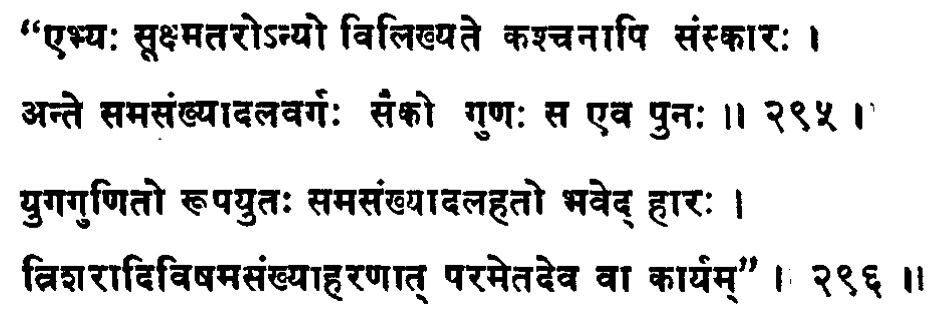}
\end{center}
Here is an English translation of the verses: 
\begin{quote}
``A subtler method, with another correction. [Retain] the first procedure involving division of four times the diameter by the odd numbers, 3, 5, etc. [But] then add or subtract it [four times the diameter] multiplied by one added to the next even number halved and squared, and divided by one added to four times the preceding multiplier [with this] multiplied by the even number halved.''
\end{quote}
Taking the diameter of a given circle as $d$, the result stated in the verses can be presented as follows:
\begin{equation}
\text{Circumference } = 4d - \frac{4d}{3}  + \frac{4d}{5}- \cdots \pm \frac{4d}{p}  \mp 4d\frac{ \left( \frac{p+1}{2} \right)^2+1}{ \left( \left[\left(\frac{p+1}{2}\right)^2 +1\right]  4 + 1 \right)\left(\frac{p+1}{2}\right)}
\end{equation}
The correction term is 
\begin{align*}
f_3(p) & = \frac{ \left( \frac{p+1}{2} \right)^2+1}{ \left( \left[\left(\frac{p+1}{2}\right)^2 +1\right]  4 + 1 \right)\left(\frac{p+1}{2}\right)}\\
&= \frac{{{p}^{2}}+2 p+5}{2 {{p}^{3}}+6 {{p}^{2}}+16 p+12}
\end{align*}
Obviously, if we let $p=2n-1$, this yields the correction term $F_3(n)$.
\subsubsection*{Remark 1}
The first correction term though not explicitly stated as such has been indicated in the rationale for the  correction term $f_2(p)$. We shall denote it by $f_1(p)$ and it is given by (see Eq.\eqref{eq20a})
$$
f_1(p)=\frac{1}{2(p+1)}.
$$
\subsubsection*{Remark 2}
As can be seen from the above description, the Kerala texts presents the correction term as a function of $p$ which is ``... the odd number at which the process ends''. In Eq.\eqref{eq1} we have presented the correction term as function of $n$ which specifies the position of the last term where the series is truncated. The relation between the two functions $f(p)$ and $F(n)$ can be stated as follows:
$$
F(n)=f(2n-1)
$$
or, as
$$
f(p) = F\left(\frac{p+1}{2}\right).
$$
\subsection{The effectiveness of the correction terms in numerical computations}

That even with small values of $n$ the third correction term $F_3(n)$ yields remarkably accurate values of $\pi$ is demonstrated in Table 1.
\renewcommand{\arraystretch}{1.5}
\begin{table}[!h]
\begin{center}
\caption{Approximate values of $\pi$ computed using truncated M\=adhava-Leibniz series with the third correction term $F_3(n)$ added (the correct digits are underlined)}
\begin{tabular}{c|c}
\hline
Number of terms used &\quad Computed value of $\pi$\\
\hline
5&\quad \underline{3.1415} 87301587302\\
6&\quad \underline{3.14159} 4274480180\\
10&\quad \underline{3.141592} 705349155\\
11&\quad \underline{3.1415926} 26657870\\
20&\quad \underline{3.14159265} 4019864\\
21&\quad \underline{3.141592653} 283544\\
30&\quad \underline{3.141592653} 615266\\
31&\quad \underline{3.1415926535} 69532\\
\hline
\end{tabular}
\end{center}
\end{table}
\renewcommand{\arraystretch}{1}

That these are indeed successively better approximations can be seen from the values given in Table 2 which are computed using the computer algebra system {\em Maxima version 5.46.0}  which is a free software released under the terms of the GNU General Public License (GPL)\footnote{To download and install the software visit the webpage {\tt https://maxima.sourceforge.io/index.html}.}. The table contains values of the errors in computing the value of $\pi$ using the truncated M\=adhava-Leibniz series and using the truncated  series with the addition of the various correction terms. These errors are denoted as follows:
\begin{align*}
E(n) & = \pi - 4\left(1 - \frac{1}{3}+\frac{1}{5} - \frac{1}{7} + \cdots + (-1)^{n-1}\frac{1}{2n-1}\right)\\
E_i(n) & = \pi - 4\left(1 - \frac{1}{3}+\frac{1}{5} - \frac{1}{7} + \cdots + (-1)^{n-1}\frac{1}{2n-1} +  (-1)^nF_i(n)\right)
\end{align*}
\section{{\em Sthaulya}: A measure of the inaccuracy of the correction term}

We have the following series for the circumference $C$ of a circle of diameter $d$:
\begin{equation}
\frac{C}{4d}=1-\frac{1}{3}+\frac{1}{5}-\cdots +(-1)^{n-1}\frac{1}{2n-1}+(-1)^n\frac{1}{2n+1}+\cdots
\end{equation}
The denominators of the terms in the series keep increasing, but extremely slowly. We will have to consider an very large number of terms even to get some reasonable accuracy. So it is obviously not convenient to obtain an accurate value of the circumference. The mathematicians of the Kerala school came up with a device to overcome this difficulty. 

Let $f(m)$, where $m$ is a positive integer, be a the function which makes the following a true statement:
\begin{equation}
\frac{\pi}{4}=1-\frac{1}{3}+ \cdots \pm \frac{1}{p}+\mp f(p)
\end{equation}
Then we must have 
\begin{align}
\frac{\pi}{4} & = 1-\frac{1}{3}+\cdots+ (-1)^{n-2}\frac{1}{2n-3}+(-1)^{n-1} f(2n-3)\label{eq10}\\
\frac{\pi}{4} & = 1-\frac{1}{3}+\cdots+(-1)^{n-2}\frac{1}{2n-3}+ (-1)^{n-1} \frac{1}{2n-1}+(-1)^{n}f(2n-1)\label{eq11}
\end{align}
From Eq.\eqref{eq10} from Eq.\eqref{eq11}, we get
\begin{equation}
(-1)^{n-1} f(2n-3) = (-1)^{n-1}\frac{1}{2n-1}+(-1)^{n}f(2n-1),
\end{equation}
Cancelling out $(-1)^{n-1}$ and rearranging the terms we get
\begin{equation}\label{eq13}
f(2n-3)+f(2n-1)-\frac{1}{2n-1} =0.
\end{equation}
So we need to find a function $f(m)$ which satisfies the recurrence relation given in Eq.\eqref{eq13}. No rational function in $n$ can satisfy this equation exactly. 

%
\renewcommand{\arraystretch}{1}
\begin{table}

\caption{Numerical values of the errors}
\medskip
\quad\hspace{-0.62in}
\begin{tabular}{|r|r|r|r|r|}
\hline
$n$ & $E(n)$ & $E_1(n)$ & $E_2(n)$ & $E_3(n)$ \\
\hline
10 & $0.0997530$ & $-2.4696534\cdot 10^{-4}$ & $2.4112190\cdot 10^{-6} $ &  $-5.1759362\cdot 10^{-8}$
\\
11 & $-0.0907231$ & $1.8593504\cdot 10^{-4}$ & $-1.5063313\cdot 10^{-6}$ & $2.6931923\cdot 10^{-8}$
\\
100 & $0.0099998$ & $-2.4996876\cdot 10^{-7}$ & $2.4990676\cdot 10^{-11}$ & $-5.7731597\cdot 10^{-15}$
\\
101 & $-0.0099007$ & $2.4261781\cdot 10^{-7}$ & $-2.3777869\cdot 10^{-11}$ & $5.3290705\cdot 10^{-15}$
\\
1000 & $9.9999975\cdot 10^{-4}$ & $-2.5000002\cdot 10^{-10}$ & $0.0$ & $0.0$
\\
1001 & $-9.9900075\cdot 10^{-4}$ & $2.4925084\cdot 10^{-10}$ & $-4.4408921\cdot 10^{-16}$ & $0.0$
\\
10000& $9.9999999\cdot 10^{-5}$ & $-2.5002223\cdot 10^{-13}$ & $0.0$ & $0.0$
\\
10001 & $-9.9990001\cdot 10^{-5} $ & $2.5002223\cdot 10^{-13}$ & $0.0$  & $0.0$\\	
\hline
\end{tabular}
\end{table}
\renewcommand{\arraystretch}{1}

In Kerala texts, the quantity
$$
f(2n-3)+f(2n-1)-\frac{1}{2n-1}
$$
is called {\em sthaulya} and it is taken as a measure of the inaccuracy of the correction term. Ideally, the {\em sthaulya} should be zero. So, the efforts were to find correction terms for which the {\em sthaulya} is as small as possible.

As we have seen in Section \ref{S4}, in Kerala texts, the correction terms are presented as functions of the number $p=2n-1$, ``... the odd number at which this process ends...'', and so, in a similar way we shall also present the {\em sthaulya} as a function of $p$ and denote it by $I(p)$. Thus
\begin{align}
I(p)& = f(2n-3)+f(2n-1)-\frac{1}{2n-1}\\
& = f(p-2)+f(p) -1/p.
\end{align}
It should be emphasised that $I(p)$ is only a measure of the inaccuracy in the correction term. It is not an estimate of the error in the value of the circumference of a circle computed using the truncated Leibniz-M\=adhava series.
\section{Kerala rationales for the correction terms}
\subsection{The first correction term}\label{Rationale}
{\em Yuktibh\=a\d sa} argues that if we choose the correction term $f(m)$ such that
\begin{equation}\label{eq14}
f(2n-3)=f(2n-1)=\frac{1}{2(2n-1)}
\end{equation}
then Eq.\eqref{eq13} is obviously satisfied, that is, the {\em sthaulya} would be exactly equal to zero. However, {\em Yuktibh\=a\d sa} immediately points out that  such a choice is absurd because it would imply that 
\begin{equation}
f(1)=f(3)=f(5)=\cdots
\end{equation}
and so $f(m)$ would have to be constant and this constant has to be equal to the fraction $$\dfrac{1}{2(2n-1)}$$ containing the variable $n$. 

{\em Yuktibh\=a\d sa} continues the argument and observes that the even integers closest to $2(2n-1)=4n-2$ are $4n-4$ and $4n$. Among these two values, as a first approximation, the even integer $4n$  is chosen and sets 
\begin{equation}\label{eq17}
f_1(2n-1)=\frac{1}{4n}
\end{equation}
or
\begin{equation}\label{eq20a}
f_1(p)=\frac{1}{2(p+1)}
\end{equation}
This is the first correction term proposed by M\=adhava. In {\em Yuktibh\=a\d sa}, the integer ``$4n$'' has been described as ``the double of the even number above the last odd number divisor'', that is, as ``$2( (2n-1) +1)$''.

\subsubsection{Computation of the {\em sthaulya}}
With this choice of $f(2n-1)$ or $f(p)$, Eq.\eqref{eq13} is  certainly not exactly satisfied.  Let us calculate how far the left side of Eq.\eqref{eq13} deviates from $0$. We calculate the {\em sthaulya} as follows. For $p\ge 3$, we have
\begin{align}
I_1(p) & = f_1(p-2) + f_1(p) - \frac{1}{p}\notag\\
     & = \frac{1}{2((p-2)+1)} + \frac{1}{2(p+1)} - 1/p\notag\\
     & = \frac{1}{2(p-1)}+\frac{1}{2(p+1)} - 1/p\notag\\
     & = \frac{1}{p^3-p} \text{  for  } p\ge 3.\label{eq19}
\end{align}
The {\em sthaulya} is inversely proportional to the third power of $p$, that is, to the third power of $n$.
\subsubsection{The choice of \texorpdfstring{$4n$}{4n} is optimum in a certain sense!}
In the above description, we observed that the two even integers closest to $2(2n-1)$ are $4n-4$ and $4n$ and we arbitrarily chose the even integer $4n$ to define the correction term. We could as well have chosen the even integer $4n-4$. We show below why this would be a bad choice. 

We consider a more general case and compute the value of $E(p)$ if we choose 
\begin{equation}\label{eq20}
f(2n-1)=\frac{1}{4n+k} = \frac{1}{2(p+1)+k}
\end{equation}
where $k$ is some constant, positive or negative. In this case we have
\begin{align}
I(p) & = \frac{1}{2(p-1)+k} + \frac{1}{2(p+1)+k} - \frac{1}{p}\notag\\
     & = -\frac{2kp+k^2-4}{p(2p+k+2)(2p+k-2)}\label{eq21}
\end{align}
This has a larger numerator than the earlier one because a term proportional to $p$ appears in the numerator. This means that, for large values of $p$, the error given by Eq,\eqref{eq21} varies inversely as $p^2$ where as the error given by Eq.\eqref{eq19} varies inversely as $p^3$. The error would vary inversely as $p^3$ only if 
$k=0$, that is only if we set
$$
f(2n-1)=\frac{1}{4n}.
$$
\subsubsection*{Remark} 
{\em Yuktibh\=a\d sa} does not discuss the general case as given Eq.\eqref{eq19}. It discusses only the special case $k=1$ and derives the expression
$$
I(p)=\frac{-2p+3}{4p^3+4p^2-3p}.
$$
\subsection{The second correction term}
We present below an outline of the rationale for the second correction term as given in Kerala texts. For easy comprehension, the outline is presented as a step-by-step process. As has been pointed out below, there are several gaps in the argument and the feeling one gets is that the  attempt to give a deductive proof to
the correction term, which had presumably been obtained inductively, was not completely successful (see \cite{Hayashi}).
\begin{enumerate}[{\bf {Step} 1.}]
\item
Firstly, we observe that the measure of inaccuracy corresponding to the first correction term $F_1(n)$ given in Eq.\eqref{eq19} is positive which means that the correction term has ``over-corrected'' the value. 
\item
To reduce the error, we must have a smaller value for the correction term and for this we must increase the denominator of the correction term.
\item
As a trial, we increase the value of the denominator by unity and consider the correction term
$$
F_{1a}(n)=\frac{1}{4n+1} =\frac{1}{2p+3}=f_{1a}(p).
$$
\item
The measure of inaccuracy is recalculated with the correction term considered in Step 3.
\begin{align}
I_{1a}
&=f_{1a}(p-2)+f_{1a}(p) - \frac{1}{p}\notag\\
&=\frac{1}{2(p-2)+3}+\frac{1}{2p+3}-\frac{1}{p}\notag\\
&=\frac{1}{2p-1}+\frac{1}{2p+3}-\frac{1}{p}\notag \\
&=\frac{2p^2+3p}{4 {{p}^{3}}+4 {{p}^{2}}-3} + \frac{2p^2-p}{4 {{p}^{3}}+4 {{p}^{2}}-3 }-\frac{4p^2+4p-3}{4 {{p}^{3}}+4 {{p}^{2}}-3 }\notag\\
&=-\frac{2 p-3}{4 {{p}^{3}}+4 {{p}^{2}}-3 p}
\end{align}
\item
This expression shows that the measure of inaccuracy has now become negative and hence it has under-corrected the value. Moreover, since the presence the term $2p$ in the numerator indicates that the error has a contribution from the value of $p$ where as in the case of $F_(n)$ there was no such contribution.
\item
The observation in Step 5 implies that we should add a quantity less than unity to the denominator of the expression for $F_1(n)$. {\em Yuktibh\=a\'sa} says: ``... a full 1 should not be added to the correction divisor''.
\item
As an improvement, the following correction term is then considered without any proper explanation (see \cite{Hayashi}):
\begin{align}
F_{1b}(n)& =\frac{1}{4n+\dfrac{1}{4n}}\\
      & = \frac{1}{(2p+2) + \frac{1}{2p+2}} =f_{1b}(p)
\end{align}
\item
The number  $4$ in the numerator in the expression $\frac{4}{4n}$ in the correction term $F_2(n)$ is explained thus: ``In this case the difference is only unity, but the difference to be obtained is four
because the numerator of the odd number is less than that of the sum of the others
by four. Therefore, both denominators of the corrections 
supposed before should be increased by the four divided by themselves.'' 

This is interpreted as follows. In the present case, that is, in the case of $f_{1b}(p)$, we have
\begin{align*}
I_{1b}
& = \Big(f_{1b}(p-2) + f_{1b}(p)\Big) -\frac{1}{p}\\
& = \left(\frac{1}{(2p-2) + \frac{1}{2p-2}}  + \frac{1}{(2p+2) + \dfrac{1}{2p+2}}\right)  - \dfrac{1}{p} \\
& \approx \left(\frac{1}{2p+\frac{1}{2p}} + \frac{1}{2p+\frac{1}{2p}} \right) -\frac{1}{p}\\
& \approx \left(\frac{4p}{4p^2+1}\right) - \frac{1}{p}\\
& \approx \left(\frac{4p^2}{p(4p^2+1)} \right)-\frac{4p^2+1}{p(4p^2+1)}.
\end{align*}
The difference between the numerators of sum of the correction terms and the numerator of the term $\frac{1}{p}$ is approximately equal to $1$. 
But in the case of $F_1(n)$, 
\begin{align*}
I_1(p) 
& = \Big(f_1(p-2) + f_1(p)\Big)-1/p\\
& = \left(\frac{1}{2p-2}+\frac{1}{2p+2}\right)-\frac{1}{p}\\
& = \left(\frac{4p^2}{p(4p^2-4)}\right) -\frac{4p^2-4}{p(4p^2 - 4)}
\end{align*}
and the difference between the numerators is $4$. Therefore, the denominators of the corrections should be increased by four divided by themselves. Thus, the next best correction term may be taken as 
$$
F_2(n)=\frac{1}{4n+\dfrac{4}{4n}}.
$$
\end{enumerate}

\subsubsection*{Remark}
To obtain a better correction term, we experimented with an expression of the form $\frac{1}{4n+k}$. From the resulting expression for $I(p)$ given in Eq.\eqref{eq21} it was  seen that the value of $k$ cannot be a constant; it must depend on $p$. Also, since 
$$
2pk=2(2n-1)k= 4nk - 2k, 
$$ 
$4n$ appears as a factor of $k$ in the numerator of $I(p)$, this factor should be cancelled out. If we choose $k$ in the form 
\begin{equation}
\frac{\star}{4n}
\end{equation}
then $4n$ is likely to cancel out from $2kp$. So we choose 
$$
k=\frac{m}{4n}=\frac{m}{2p+2}.
$$
With this value for $k$, the correction term becomes
$$
f(p)=\frac{1}{ (2p+2) + \dfrac{m}{2p+2}}
$$
The error estimate can be shown to be
$$
I(p) = - \frac{\left( 4 m-16\right)  {{p}^{2}}+{{m}^{2}}+8 m+16}{16 {{p}^{5}}+\left( 8 m-32\right)  {{p}^{3}}+\left( {{m}^{2}}+8 m+16\right)  p}.
$$
The numerator is proportional to $p^2$ unless $4m-16=0$ or $m=4$. Thus to optimize the error we have to choose $m = 4$. The correction term now takes the form
\begin{equation}\label{eq24a}
f_2(p) = \frac{1}{ (2p+2) + \dfrac{4}{2p+2}}
\end{equation}
or
\begin{align*}
f_2(2n-1)&=\frac{1}{4n+\dfrac{4}{4n}}\\
&=\frac{n}{4n^2+1}
\end{align*}
which is precisely M\=adhava's second correction term.

\subsubsection{Measure of inaccuracy}
The measure of inaccuracy corresponding to the correction term $f_2(p)$ is given by
\begin{align}
I_2(p)& = f_2(p-2) + f_2(p) - \frac{1}{p}\notag\\
    & = \frac{1}{ (2(p-2)+2) + \dfrac{4}{2(p-2)+2}} +
        \frac{1}{(2p+2) + \dfrac{4}{2p+2}}- \frac{1}{p}\notag\\
    & = -\frac{4}{p^5+4p}\text{ for } p\ge 3\label{eq24}
\end{align}
This is inversely proportional to the fifth power of $p$.
\subsubsection{This {\em sthaulya} is also optimum in the same sense!}
We show that the error cannot be reduced by choosing the correction term in the following form with an integer constant $k$:
\begin{equation}\label{eq25}
f(2n-1)=\frac{1}{4n + \dfrac{4}{4n+k^\prime}},
\end{equation}
or
$$
f(p)=\frac{1}{(2p+2)+\dfrac{4}{(2p+2)+k^\prime}},
$$
In this case the measure of inaccuracy would be 
\begin{align}
I(p)& = f(p-2) + f(p) -1/p \notag \\
& = \frac{2 k^\prime p+{{k^\prime}^{2}}-16}{4 {{p}^{5}}+4 k^\prime {{p}^{4}}+{{k^\prime}^{2}} {{p}^{3}}+\left( 16-{{k^\prime}^{2}}\right)  p}\label{eq27}
\end{align}
Because of the presence of the term $2k^\prime p$, this has a larger numerator than the numerator of the error expression in Eq.\eqref{eq24}. It is minimum when $k^\prime=0$. 
\subsection{The third correction term}
As has already been pointed out, the texts of Kerala school of mathematics do not give any rationale for the third correction term. It appears that there is simple and straight forward rationale for this choice which so far has not been noticed by experts who have written extensively on the correction terms. In the numerator of the expression in Eq.\eqref{eq27} there is the term $2k^\prime p=2k^\prime(2n-1)= 4nk^\prime$. To cancel the factor $4n$ we should choose $k^\prime$ also in the same form as as the form for $k$, namely 
$$
k^\prime = \frac{m}{4n}.
$$
To get a third correction we substitute this value of $k^\prime$ in Eq.\eqref{eq25} and take 
$$
f(p)=\frac{1}{ (2p+2) + \dfrac{4}{ (2p+2) + \dfrac{m}{(2p+2)}}}. 
$$
This correction term can be shown to produce the following value for the {\em sthaulya}:
$$
I(p)= \frac{\left( 4 m-64\right)  {{p}^{2}}+{{m}^{2}}+16 m+64}{16 {{p}^{7}}+\left( 8 m-16\right)  {{p}^{5}}+\left( {{m}^{2}}+8 m+64\right)  {{p}^{3}}+\left( -{{m}^{2}}-16 m-64\right)  p}.
$$
Here also $p^2$ appears in the numerator unless $m=16$. Thus to obtain an optimum error estimate we choose $m=16$ and get the following correction term:
\begin{align}
f_3(p) &= \frac{1}{ (2p+2) + \dfrac{4}{ (2p+2) + \dfrac{16}{(2p+2)}}}\label{28a}\\
&=\frac{ \left( \dfrac{p+1}{2} \right)^2+1}{ \left( \left[\left(\dfrac{p+1}{2}\right)^2 +1\right]  4 + 1 \right)\left(\dfrac{p+1}{2}\right)}\notag
\end{align}
or
\begin{align}
f_3(2n-1) & = \frac{1}{ 4n + \dfrac{4}{ 4n + \dfrac{16}{4n}}}\\
& = \frac{n^2+1}{n(4n^2+5)}
\end{align}
which is precisely M\=adhava's third correction term.
\subsubsection{Measure of inaccuracy}
The value of the measure of inaccuracy is
\begin{align}
I_3(p) & = f_3(p-2) + f_3(p) -1/p\notag\\
& = \frac{36}{{{p}^{7}}+7 {{p}^{5}}+28 {{p}^{3}}-36 p}\,\,\text{for } p\ge 3\label{eq30a}
\end{align}
This is inversely proportional to the seventh power of $p$.
\section{More on the measure of inaccuracy}
\subsection{{\em Sthaulya} as an upper bound of the error in the computed value of \texorpdfstring{$\frac{\pi}{4}$}{pi/4}}
Letting $f(p)$ denote one of the correction terms $f_1(p)$ or $f_2(p)$ or $f_3(p)$ and setting $f(-1)=0$, we consider the following series:
\begin{equation}\label{eq32}
-\bigg(f(-1) + f(1) - 1\bigg) + \left(f(1)+f(3) -  \frac{1}{3}\right) - \left(f(3) + f(5) -\frac{1}{5}\right) +\cdots
\end{equation}
that is, the series
\begin{equation}\label{eq31a}
-I(1) +I(3) - I(5) + I(7) - \cdots
\end{equation}
where $I(p)$ is the associated {\em sthaulya}, that is, the measure of inaccuracy, given by
$$
I(p)=f(p-2) + f(p)-\frac{1}{p}.
$$
The $n$-th  term of the series is 
$$
I(2n-1) = f(2n-3)+f(2n-1) -\frac{1}{2n-1}.
$$
From Eqs.\eqref{eq19}, \eqref{eq24} and \eqref{eq30a}, we see that, for all odd $p\ge 3$, we have:
\begin{itemize}
\item
$I_1(p) = \dfrac{1}{p^3-p} > 0$
\item
$I_2(p) = - \dfrac{4}{p^5+4p} <0$
\item
$I_3(p) = \dfrac{36}{p^7+7p^5+28p^3-36p}>0$
\end{itemize}
Also in all the three cases we have
$$
\lim_{p\rightarrow \infty} I(p) =0
$$
and 
$$
|I(p+2)| < |I(p)|.
$$
We conclude that the series in Eq.\eqref{eq31a} is an alternating series whose terms monotonically decrease to zero. 

The $n$-th partial sum of the series Eq.\eqref{eq32} is
\begin{align}
s_n
& = -\bigg(f(-1) + f(1) - 1\bigg) + \left(f(1)+f(3) -  \frac{1}{3}\right) - \left(f(3) + f(5) -\frac{1}{5}\right) +\cdots\notag\\
& \qquad(-1)^{n-1}\left( f(2n-3)+ f(2n-1)-\frac{1}{2n-1}\right)\notag\\
& = 1-\frac{1}{3}-\frac{1}{5}- \cdots + (-1)^{n-1}\frac{1}{2n-1} +(-1)^nf(2n-1)
\end{align}
Letting $n$ tend to $\infty$, we have
\begin{align}
\lim_{n\rightarrow \infty} s_n 
& = 1 -\frac{1}{3} + \frac{1}{5} - \cdots (-1)^{n-1} \frac{1}{2n-1}+ \cdots\notag\\
& = \frac{\pi}{4}\,\, \text{ by Eq.\eqref{eq1a}}\label{eq36a}
\end{align}
From properties of alternating series we have
$$
\left|\frac{\pi}{4}-s_n\right| \le |I(2n+1)|
$$
that is,
$$
\left|\frac{\pi}{4}-s_n\right| \le \left|f(2n-1) +f(2n+1) - \frac{1}{2n+1}\right|
$$
\subsubsection*{Another error estimate using {\em sthaulya}}
The sequence $\{I(2n+1) - I(2n-1): n=1,2,\cdots\}$  monotonically converge to zero. Hence (see \cite{Calabrese}) we have the following bounds for the error in the computation of $\frac{\pi}{4}$.
$$
\frac{1}{2}|I(2n+1)|<\left|\frac{\pi}{4} - s_n\right|<\frac{1}{2}|I(2n-1)|.
$$
\subsection{The correction terms as convergents of a continued fraction}
Let $p=2n-1$ and the remainder after $n$ terms in the Leibniz series be $C(p)$ so that we may write
$$
\frac{\pi}{4}=1-\frac{1}{3}+\frac{1}{5}-\cdots\pm \frac{1}{p}\mp C(p).
$$
Then the remainder term has the following continued fraction expansion (see \cite{Dukta}, \cite{Rajagopal}) :
\begin{equation}
C(p)= \frac{1}{(2p+2) \,\, +}\,\, \frac{2^2}{(2p+2)\,\, +}\,\, \frac{4^2}{(2p+2)\,\, +}\,\, \frac{6^2}{(2p+2)\,\, +}\,\,\frac{8^2}{(2p+2)\,\,+}\cdots
\end{equation}
The continued fraction $C(p)$ can be put in different forms. For example we have:
$$ 
C(p)=\frac{1}{2}\left( \frac{1}{(p+1)\,\, +}\,\, \frac{1^2}{(p+1)\,\,+}\,\,\frac{2^2}{(p+1)\,\,+}\,\,\frac{3^2}{(p+1)\,\,+}\,\,\frac{4^2}{(p+1)\,\,+}\cdots\right)
$$

From Eqs.\eqref{eq20a}, \eqref{eq24a} and \eqref{28a} we can see that the the correction terms $f_1(p)$, $f_2(p)$ and $f_3(p)$ are the first three convergents of the infinite continued fraction $C(p)$.
\subsection{Extensions of M\=adhava's results: The fourth and fifth correction terms}
Extending  the rationale of the mathematicians of the Kerala school, we can derive further better correction terms. The fourth and the fifth correction terms with the associated {\em sthaulya}'s are given below:
\begin{align*}
f_4(p)&= \frac{1}{2 p+2+\dfrac{4}{2 p+2+\dfrac{16}{2 p+2+\dfrac{36}{2 p+2}}}}\\
&= \frac{{{p}^{3}}+3 {{p}^{2}}+16 p+14}{2 {{p}^{4}}+8 {{p}^{3}}+40 {{p}^{2}}+64 p+48}\\
I_4(p)& = f_4(p-2)+f_4(p)-\frac{1}{p}\\
& = - \frac{576}{{{p}^{9}}+24 {{p}^{7}}+192 {{p}^{5}}-64 {{p}^{3}}+576 p}\\
& = -\frac{576}{p\, \left( {{p}^{4}}-4 {{p}^{3}}+20 {{p}^{2}}-32 p+24\right) \, \left( {{p}^{4}}+4 {{p}^{3}}+20 {{p}^{2}}+32 p+24\right) }\\
f_5(p)&=\frac{1}{2 p+2+\dfrac{4}{2 p+2+\dfrac{16}{2 p+2+\dfrac{36}{2 p+2+\dfrac{64}{2 p+2}}}}}\\
&=\frac{{{p}^{4}}+4 {{p}^{3}}+35 {{p}^{2}}+62 p+94}{2 {{p}^{5}}+10 {{p}^{4}}+80 {{p}^{3}}+200 {{p}^{2}}+368 p+240}\\
I_5(p)&=f_5(p-2)+f_5(p)-\frac{1}{p}\\
&=\frac{14400}{{{p}^{11}}+55 {{p}^{9}}+968 {{p}^{7}}+3520 {{p}^{5}}+9856 {{p}^{3}}-14400 p}\\
& = \frac{14400}{(p^3-p) \, \left( {{p}^{4}}-4 {{p}^{3}}+36 {{p}^{2}}-64 p+120\right) \, \left( {{p}^{4}}+4 {{p}^{3}}+36 {{p}^{2}}+64 p+120\right) }\\
\end{align*}
\subsection{Some elementary properties of the measures of inaccuracy}
Let Let the reduced form of the rational function $I_k(p)$ be 
$I_k=\dfrac{N_k(p)}{D_k(p)}.$
The following table shows the factorisation of the denominator polynomial $D_k(p)$ for the first few values of $k$.
\renewcommand{\arraystretch}{1.4}
\begin{table}
\begin{center}
\caption{Factorisations of the denominators of the measures of inaccuracy}
\begin{tabular}{c|l}
\hline
$k$ & $N_k(p)$\\
\hline
$1$ & $p(p-1)(p+1)$\\
$2$ & $p\, \left( {{p}^{2}}-2 p+2\right) \, \left( {{p}^{2}}+2 p+2\right)$\\
$3$ & $\left( p-1\right)  p\, \left( p+1\right) \, \left( {{p}^{2}}-2 p+6\right) \, \left( {{p}^{2}}+2 p+6\right)$\\
$4$ & $ p\, \left( {{p}^{4}}-4 {{p}^{3}}+20 {{p}^{2}}-32 p+24\right) \times $ \\
 &\,\qquad $ \left( {{p}^{4}}+4 {{p}^{3}}+20 {{p}^{2}}+32 p+24\right)$\\
$5$ & $\left( p-1\right)  p\, \left( p+1\right) \, \left( {{p}^{4}}-4 {{p}^{3}}+36 {{p}^{2}}-64 p+120\right) \times $\\
    &\,\qquad$ \left( {{p}^{4}}+4 {{p}^{3}}+36 {{p}^{2}}+64 p+120\right) $\\
$6$ & $p\, \left( {{p}^{2}}-2 p+10\right) \, \left( {{p}^{2}}+2 p+10\right) \, \left( {{p}^{4}}-4 {{p}^{3}}+52 {{p}^{2}}-96 p+72\right) \times $\\
  & \,\qquad$ \left( {{p}^{4}}+4 {{p}^{3}}+52 {{p}^{2}}+96 p+72\right)$\\
$7$ & $\left( p-1\right)  p\, \left( p+1\right) \, \left( {{p}^{6}}-6 {{p}^{5}}+106 {{p}^{4}}-384 {{p}^{3}}+2080 {{p}^{2}}-3408 p+5040\right)\times $\\
 & \,\qquad$ \left( {{p}^{6}}+6 {{p}^{5}}+106 {{p}^{4}}+384 {{p}^{3}}+2080 {{p}^{2}}+3408 p+5040\right)$\\
\hline
\end{tabular}
\end{center}
\end{table}
From the expressions for the measures of inaccuracy considered so far, we may make some elementary observations about these expressions.
\begin{enumerate}
\item
$N_k(p)$ is a constant and is equal to $(-1)^{k+1} (k!)^2$.
\item
$D_k(p)$ is a monic polynomial of degree $2k-1$. 
\item
$D_k(p)$ contains only odd powers of $p$.
\item
$D_k(0)=0$.
\item
The coefficient of $p^{2k-3}$ in $D_k(p)$  is 
$$
\frac{1}{3}(k-2)k(2k+1)= 2(1^2+2^2+\cdots+k^2) - (1+2+3+\cdots+2k).
$$
\item
The coefficient of $p$ in $D_k(p)$ is $(-1)^k (k!)^2$.
\item
If $k$ is an odd integer then $(p-1)p(p+1)$ is a factor of $D_k(p)$. Also $D_k(p)$ can be factorised in the form $(p-1)p(p+1)\psi(p)\psi(-p)$ where $\psi(p)$ is a monic polynomial of degree $k-1$ in $p$
\item
If $k$ is an even integer, $D_k(p)$ can be factorised in the form $p\psi(p)\psi(-p)$ where $\psi(p)$ is monic polynomial of degree $k$ in $p$. 
\end{enumerate}
\section{More fastly converging series for \texorpdfstring{$\pi$}{pi}}
The expressions for the measures of inaccuracy have been put to use by mathematicians of the Kerala school to obtain convergent series for $\pi$ which converge more fastly than the Madhava-Leibniz series. The starting point for obtaining such series is the following series representation of $\pi$ given in Eq.\eqref{eq36a}:
$$
\frac{\pi}{4}= \bigg(1-f(1)\bigg) + \sum_{n=2}^\infty (-1)^{n}\left(f(2n-3)+f(2n-1)-\frac{1}{2n-1}\right)
$$
where $f(m)$ is a measure of inaccuracy.  Now we consider the first three measures of inaccuracies. 

For more details on the various series representations of $\pi$ considered in this section, the reader may consult \cite{Rajagopal}, \cite{Raju} pp.176-178 and \cite{Yukti} pp.205-206.
\subsection{Series using Madhava's correction terms}
\subsubsection*{Case $k=1$:}
\begin{align}
\frac{\pi}{4}
& = \bigg(1-f_1(1)\bigg) + \sum_{n=2}^\infty (-1)^{n}\left(f_1(2n-3)+f_1(2n-1)-\frac{1}{2n-1}\right) \notag\\
& = \bigg(1-\frac{1}{4}\bigg) + \sum_{n=2}^\infty(-1)^n\left( \frac{1}{4n-4}+\frac{1}{4n}-\frac{1}{2n-1}\right)\notag\\
& = \frac{3}{4} + \sum_{n=2}^\infty (-1)^n\left(\frac{1}{(2n-1)^3 -(2n-1)}\right)\notag\\
& = \frac{3}{4} + \frac{1}{3^3-3}-\frac{1}{5^3-5}+\frac{1}{7^3-7} - \cdots\label{eq38}
\end{align}
The verse stating this result is as follows (\cite{Tantra}, Chapter II, verse 290):
\begin{center}
\includegraphics[width=8cm]{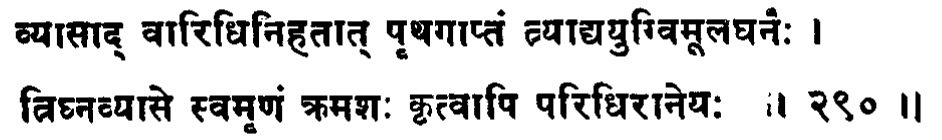}
\end{center}
The English translation of the verse is given below:
\begin{quote}
``Four times the diameter is divided by the cubes of [odd numbers] 3, etc., minus
the numbers [lit. roots], to obtain separate quotients. To thrice the diameter,
alternately add and subtract [the quotients], to obtain the circumference.''
\end{quote}
\subsubsection*{Case $k=2$:}
\begin{align}
\frac{\pi}{4}
& = \bigg(1-f_2(1)\bigg) + \sum_{n=2}^\infty (-1)^{n}\left(f_2(2n-3)+f_2(2n-1)-\frac{1}{2n-1}\right) \notag\\
& = \bigg(1- \frac{1}{5}\bigg) + \sum_{n=2}^\infty (-1)^n \left( \frac{n-1}{4(n-1)^2+1} + \frac{n}{4n^2+1} - \frac{1}{2n-1}\right)\notag\\
& = \frac{4}{5} + \sum_{n=2}^\infty (-1)^{n+1}\frac{4}{(2n-1)^5 + 4(2n-1)}\notag\\
& = \frac{4}{1^5 + 4\cdot 1} - \frac{4}{3^5+4\cdot 3}+\frac{4}{5^5+4\cdot 5} - \frac{4}{7^5+4\cdot 7}+ \cdots\label{eq39}
\end{align}
The verse stating this result is as follows (\cite{Tantra}, Chapter II, verses 287-188):
\begin{center}
\includegraphics[width=8cm]{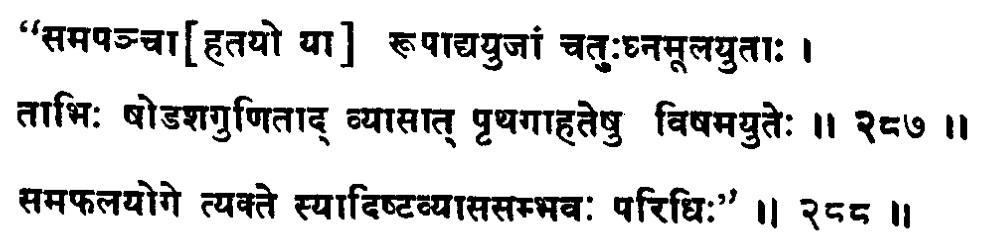}
\end{center}
The English translation of the verse is given below:
\begin{quote}
``The fifth powers of 1, etc., plus four times the number; with that, divide 16
times the diameter separately for successive odd numbers and alternately add
and subtract. The circumference is obtained for the desired diameter.''
\end{quote}
\subsubsection*{Case $k=3$:}

\noindent This case is not discussed in Kerala texts.
\begin{align}
\frac{\pi}{4}
& = \bigg(1-f_3(1)\bigg) + \sum_{n=2}^\infty (-1)^{n}\left(f_3(2n-3)+f_3(2n-1)-\frac{1}{2n-1}\right) \notag\\
& = \bigg(1-\frac{2}{9}\bigg) + \sum_{n=2}^\infty (-1)^n\left( \frac{(n-1)^2+1}{(n-1)(4(n-1)^2+5)} + \frac{n^2+1}{n(4n^2+5)}-\frac{1}{2n-1}\right)\notag\\
& = \frac{7}{9} + \sum_{n=2}^\infty (-1)^n\left( \frac{36}{(2n-2)(2n-1)(2n)(((2n-1)^2+4)^2+20)}\right)\notag\\
& = \frac{7}{9} 
+ \frac{36}{2\cdot 3 \cdot 4 ((3^2+4)^2+20)}
- \frac{36}{4\cdot 5 \cdot 6 ((5^2+4)^2+20)}
 \cdots\label{eq40}
\end{align}
The series can also be put in the following form (see \cite{Srinivas}):
$$
\frac{\pi}{4} = \frac{7}{9} 
+ \frac{36}{(3^3-3)(2^2+5)(4^2+5)} 
-\frac{36}{(5^3-5)(4^2+5)(6^2+5)}
+\cdots
$$
\subsubsection*{Remark}
Figure 1 shows the common logarithms of the errors in the calculated value of $\frac{\pi}{4}$ computed using the Madhava series (Eq.\eqref{eq1a}) and the  transformed series derived using the three correction terms (Eqs.\eqref{eq38}, \eqref{eq39} and \eqref{eq40}).
\begin{figure}
\begin{center}
\includegraphics[width=10cm]{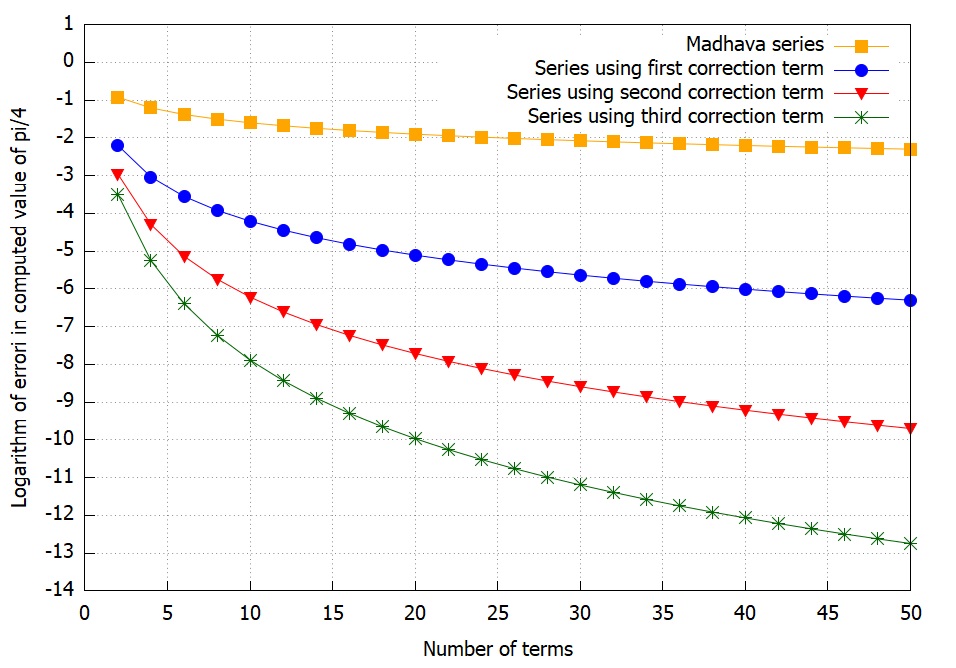}\\[2mm]
\caption{Number of terms vs log-error}
\end{center}
\end{figure}
\subsection{Series based on  ``nonoptimal'' correction term} 
We have seen than the correction terms discovered by Madhava are optimal in a certain sense. The nonoptimal correction terms have also been used by Kerala mathematicians to obtain infinite series representations for $\pi$. In fact, any function $f(p)$ which make the series in Eq.\eqref{eq32} an alternating series whose terms monotonically converge to zero can be used to construct an infinite series representation of $\pi$. It appears that Kerala mathematicians were aware of this general result.  For example, if we choose the following nonoptimal correction term
$$
f(p)=\frac{1}{2p}
$$
then we have 
$$
f(p-2)+f(p)-\frac{1}{p}= \frac{1}{(p-1)^2 -1 }
$$
With this choice of $f(p)$ the series in Eq.\eqref{eq32} can be shown to be an alternating series whose terms monotonically converge to zero. Hence we have:
\begin{align}
\frac{\pi}{4} 
& = -(\frac{1}{2} -1 ) + \sum_{n=1}^\infty (-1)^n\left( \frac{1}{(2n)^2-1}\right)\notag \\
& = \frac{1}{2} +\frac{1}{2^2-1} -\frac{1}{4^2-1}+\frac{1}{6^2-1} - \cdots\label{eq41}
\end{align}
The verse stating this result is as follows (\cite{Tantra}, Chapter II, verse 292):
\begin{center}
\includegraphics[width=8cm]{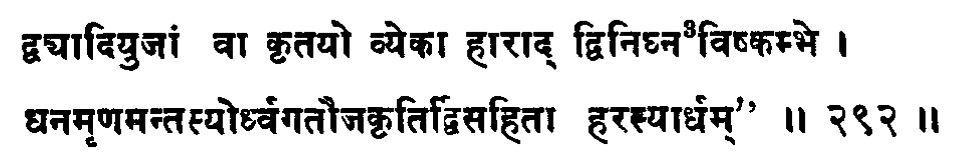}
\end{center}
English translation (see  \cite{Rajagopal2}):
\begin{quote}
``Four times the diameter is first divided successively
by the even numbers beginning with 2 squared and diminished by unity. Next,
four times the diameter is divided by double the result of squaring the odd
number next greater than the last even number and adding 2 to it. Lastly, when
all the primary quotients are alternately added to and subtracted from twice
the diameter, then results an approximation to the circumference, with a correction which is the last term in the series thus defined.''
\end{quote}
\subsection{More series representations of \texorpdfstring{$\pi$}{pi}}
{\em Yuktibhasha} (\cite{Yukti}, p.206) quotes the following two more series representations of $\pi$ as given in {\em Yuktidipika} without giving any rationale for the truth of these representations. 
\begin{align}
\frac{\pi}{8}
& = \frac{1}{2^2-1}+\frac{1}{6^2-1}+\frac{1}{10^2-1}+\cdots\label{eq42}\\
\frac{\pi}{8}
& = \frac{1}{2}- \frac{1}{4^2-1}-\frac{1}{8^2-1}-\frac{1}{12^2-1}-\cdots
\end{align}
The first of these can be obtained thus.
\begin{align*}
\frac{\pi}{4}
& = 1 - \frac{1}{3}+\frac{1}{5}-\frac{1}{7}+\frac{1}{9} -\frac{1}{11}+\cdots \\
& = \left(1 - \frac{1}{3}\right) + \left(\frac{1}{5}-\frac{1}{7}\right) + \left(\frac{1}{9} -\frac{1}{11}\right) +\cdots \\
& = \frac{2}{1\cdot 3} + \frac{2}{5\cdot 7} + \frac{2}{9\cdot 11} +\cdots \\
& = 2\left( \frac{1}{2^2-1} + \frac{1}{6^2-1} +\frac{1}{10^2-1}+\cdots \right)\\
\frac{\pi}{8} 
& = \frac{1}{2^2-1} + \frac{1}{6^2-1} +\frac{1}{10^2-1}+\cdots
\end{align*}
The latter is obtained by subtracting Eq.\eqref{eq42} from Eq.\eqref{eq41}.
The verses stating this result is as follows (\cite{Tantra}, Chapter II, verses 293-294):
\begin{center}
\includegraphics[width=8cm]{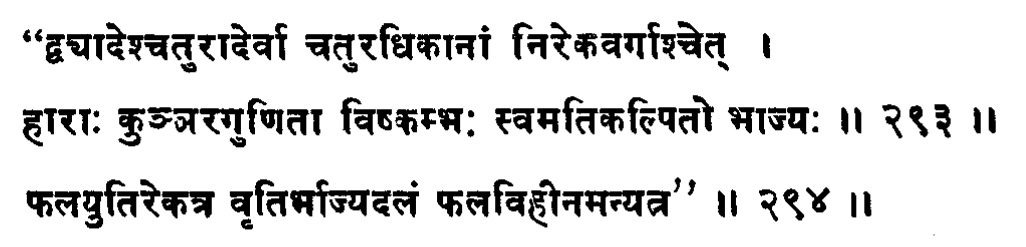}\\
\end{center}
\subsection{A series representation of \texorpdfstring{$\pi$}{pi} due to Putumana Somay\=ji}
Puthumana Somayaji (c.1660–1740), a 17th-century astronomer-mathematician from Kerala, perhaps continuing the work of M\=adhava and others, discovered a new series representation for $\pi$. It has been included in his seminal work titled {\em Kara\d napaddhati}. Venkiteswara Pai and others who have brought out an English translation with detailed notes of {\em Kara\d napaddhati} have appropriately named this series as the ``Putumana Somay\=aj\=i series'' (see \cite{Pai} p.151).  

The Putumana Somay\=aj\=i series is encoded in the following Sanskrit verse.
\begin{center}
\includegraphics[width=8cm]{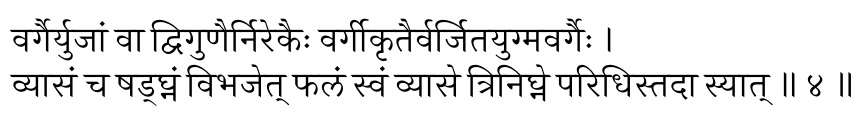}
\end{center}
Here is an English translation of the verse.
\begin{quote}
``Or, from the square of even numbers multiplied by two, subtract one, and from
the square [of that] subtract the square of the same (even number). Divide the
diameter multiplied by six by the above [quantities]. When [the sum of] these is
added to three times the diameter, the result will be the circumference.''
\end{quote}
In modern notations, this result is equivalent to the following expression for $\pi$.
$$
\frac{\pi}{4} = \frac{3}{4}+\frac{3}{2}\left( 
\frac{1}{(2\cdot 2^2 -1)^2 - 2^2}+\frac{1}{(2\cdot 4^2-1)^2-4^2} + \frac{1}{(2\cdot 6^2 - 1)^2 - 6^2}+\cdots\right)
$$
That this result can be derived from Eq.\eqref{eq38} as follows. From Eq.\eqref{eq38} we have 
\begin{align*}
\frac{\pi}{4} & = \frac{3}{4} + \frac{1}{3^3-3}-\frac{1}{5^3-5} + \frac{1}{7^3-7} - \frac{1}{9^3-9}+\cdots\\
& = \frac{3}{4} +\sum_{n=1}^\infty \left ( \frac{1}{(4n-1)^2-(4n-1)} - \frac{1}{(4n+1)^3-(4n+1)}\right)
\end{align*}
With some algebraic manipulations it can be shown that 
$$
\frac{1}{(4n-1)^2-(4n-1)} - \frac{1}{(4n+1)^3-(4n+1)} =\frac{3}{2}\times \frac{1}{ (2(2n)^2-1)^2 - (2n)^2}
$$
from which we get Putumana Somay\=aj\=i series.
\section{Hayashi's conjecture on how the three correction terms were discovered}\label{Hayashi}
In a paper published in 1990 (\cite{George} pp.132-133, \cite{Hayashi}) a group of three Japanese scholars T. Hayashi, T. Kusuba and M. Yanoa,  proposed an ingenious method by which M\=adhava might have obtained the three correction terms. Their proposal was based on two assumptions: 
\begin{itemize}
\item Madhava used $\frac{355}{113}$ as the value of $\pi$.
\item
Madhava used euclidean algorithm for division.
\end{itemize}
Their reasoning proceeds as follows. Let us write
$$
S(n) = \left| 1 - \frac{1}{3} +\frac{1}{5}-\frac{1}{7}+ \cdots +\frac{(-1)^{n-1}}{2n-1} - \frac{\pi}{4}\right|
$$
and compute the following values assuming that $\pi$ is $\frac{355}{113}$.
\begin{align*}
S(1) & = 97/452\\
S(2) & = 161/1356\\
S(3) & = 551/6780\\
S(4) & = 2923/47460\\
S(5) & = 21153/427140\\
\end{align*}
\subsubsection*{First correction term}
The values of $S(1), S(2), \cdots$ are expressed as fractions with $1$ as numerator as follows:
\begin{align*}
S(1) & = \frac{1}{4 + \frac{64}{97}} \\
S(2) & = \frac{1}{8 + \frac{68}{161}} \\
S(3) & = \frac{1}{12 +\frac{168}{551}} \\
S(4) & = \frac{1}{16 +\frac{ 692}{2923}} \\
S(5) & = \frac{1}{20 +\frac{ 4080}{21153}} 
\end{align*}
Ignoring the fractional parts in the denominator, the following approximations are obtained.
\begin{align*}
S(1) &  \approx \frac{1}{4}\\
S(2) &  \approx \frac{1}{8}\\
S(3) &  \approx \frac{1}{12}\\
S(4) &  \approx \frac{1}{16}\\
S(5) &  \approx \frac{1}{20}
\end{align*}
This suggests the following first approximation to $S(n)$ which is the first correction term $F_1(n)$. 
$$
S(n) \approx \frac{1}{4n}.
$$
\subsubsection*{Second correction term}
The fractions that were ignored are now expressed with $1$ as numerator. The fractional parts in the denominators are ignored to get the following approximations:
\begin{align*}
\frac{64}{97} & = \frac{1}{ 1 + \frac{33}{64} } \approx \frac{1}{1}\\
\frac{68}{161} & = \frac{1}{ 2 + \frac{25}{68} } \approx \frac{1}{2}\\
\frac{168}{551} & = \frac{1}{ 3 + \frac{47}{168} } \approx \frac{1}{3}\\
\frac{ 692}{2923} & = \frac{1}{ 4 + \frac{155}{692} } \approx \frac{1}{4}\\
\frac{ 4080}{21153} & = \frac{1}{ 5 + \frac{753}{4080} } \approx \frac{1}{5}
\end{align*}
This yields the next better approximation to $S(n)$ which is exactly the same as the second correction term $F_2(n)$.
$$
S(n)\approx \frac{1}{4n+\frac{1}{n}} = \frac{n}{4n^2+1}
$$
\subsubsection*{Third correction term}
The fractions that are ignored are once again expressed as fractions with $1$ as numerator to get the following approximations.
\begin{align*}
\frac{33}{64} & = \frac{1}{1+ \frac{31}{33}} \approx \frac{1}{1} \\
\frac{25}{68} & = \frac{1}{2+ \frac{18}{25}} \approx \frac{1}{2} \\
\frac{47}{168} & = \frac{1}{3+ \frac{27}{47}} \approx \frac{1}{3} \\
\frac{155}{692} & = \frac{1}{4+ \frac{72}{155}} \approx \frac{1}{4} \\
\frac{753}{4080} & = \frac{1}{5+ \frac{315}{753}} \approx \frac{1}{5} 
\end{align*}
Finally this gives a still better approximation to $S(n)$ which is the third correction term $F_3(n)$ attributed to Madhava.
$$
S(n) \approx \frac{1}{4n +\dfrac{1}{n+\dfrac{1}{n}}} = \frac{n^2+1}{n(4n^2+5)}
$$
\subsubsection*{Remark}
The above described derivations of M\=adhava's correction terms is based explicitly on the assumption that M\=adhava knew the approximate value $\frac{355}{113}$ of $\pi$. Hayashi and others have not cited any documentary evidence in support of this assumption except a claim that this approximate value of $\pi$ was known in India at least as early as the ninth century. However, Hayashi's derivation of the correction terms can still be shown to be valid even if we use the following approximation of $\pi$ (see \cite{Gupta}, \cite{MacTutor}).
$$
\pi = \frac{62832}{20000}.
$$
This approximate value of $\pi$ should definitely have been known to M\=adhava as it has been stated in \=Arybha\d t\=iya and \=Aryabha\d t\=iya has been extensively studied and commented upon by M\=adhava's student Va\d ta\'{ss}eri Parame\'{s}varan Namp\=utiri.

With this value of $\pi$ we have %
\begin{align*}
S(1) & = 1073/5000\\
S(2) & = 1781/15000\\
S(3) & = 1219/15000\\
S(4) & = 6467/105000\\
S(5) & = 15599/315000
\end{align*}
We also have 
\begin{align*}
S(1) & = \frac{1}{4  + \dfrac{1}{1 + \dfrac{1}{1 + \frac{343}{365}}}}\\
S(2) & = \frac{1}{8  + \dfrac{1}{2 + \dfrac{1}{2 + \frac{198}{277}}}}\\
S(3) & = \frac{1}{12 + \dfrac{1}{3 + \dfrac{1}{3 + \frac{63}{103}}}}\\
S(4) & = \frac{1}{16 + \dfrac{1}{4 + \dfrac{1}{4 + \frac{108}{355}}}}\\
S(5) & = \frac{1}{20 + \dfrac{1}{5 + \dfrac{1}{6 + \frac{26}{499}}}}
\end{align*}
Ignoring the fractions and the ``anomalous'' number $6$ in $S(5)$, these continued fractions suggests the following approximate expression for $S(n)$
$$
S(n)\approx \frac{1}{4n + \dfrac{1}{n + \dfrac{1}{n}}}
$$
and the following successively better approximations to $S(n)$:
\begin{alignat*}{2}
F_1(n) & = \frac{1}{4n} &&\\
F_2(n) & = \frac{1}{4n +\dfrac{1}{n}} &&= \frac{n}{4n^2+1}\\
F_3(n) & = \frac{1}{4n + \dfrac{1}{n + \dfrac{1}{n}}} &&= \frac{n^2+1}{4n^3+5n}
\end{alignat*}
\subsubsection*{Derivation of the correction terms suggested by Yushkevich}
Adolph P. Yushkevich (also spelt as Youschkevitch  and Juschkewitsch), a Soviet historian of mathematics, suggested a possible way in which Madhava might have arrived at his correction terms. But his suggested derivation is unsatisfactory on two counts: firstly, he assumes that M\=adhava used $\frac{377}{120}$ as the exact value of $\pi$ and, secondly, it looks quite artificial. After suggesting that M\=adhava might have obtained $F_1(n)$ and $F_2(n)$ more or less as conjectured by Hayashi, Yushkevich suggests that M\=adhava derived $F_3(n)$ by the following argument.
Multiply the numerator and denominator of $F_2(n)$ by $n$ to get 
$$
F_2(n) = \frac{n^2}{n(4n^2+1)}.
$$
Then $F_3(n)$ was obtained by adding the numerators and denominators of $F_1(n)$ and $F_2(n)$ separately to get
\begin{align*}
F_3(n) &=\frac{\text{[Numerator of $F_1(n)$] $+$ [Numerator of $F_2(n)$]}}{\text{[Denominator of $F_1(n)$] $+$ [Denominator of $F_2(n)$]}}\\
& = \frac{1+n^2}{4n+n(4n^2+1)}\\
& = \frac{n^2+1}{n(4n^2+5)}
\end{align*}
That indeed is really an artificial derivation of $F_3(n)$! (For more details on Yushkevic's derivation see \cite{Hayashi}.)
\section{Were all these results actually used to compute \texorpdfstring{$\pi$}{pi}?}
We conclude this paper by raising the question whether M\=adhava and his disciples did actually use the various infinite series representations of $\pi$ and the various correction terms discussed earlier in this paper to compute the value of $\pi$. The author conjectures that most probably none of them were actually used to compute the value of $\pi$! Instead, M\=adhava is likely to have used the following rapidly convergent series also known to him (see \cite{Gupta2}):
\begin{equation}\label{SankaraVarman}
\pi=\sqrt{12}\left(1 - \frac{1}{3\cdot 3} + \frac{1}{5\cdot 3^2} -\frac{1}{7\cdot 3^3}+\cdots\right)
\end{equation}
In support of this conjecture it may be pointed out that {\em Sadratnamala}, a work of Sankara Varman (1774–1839) who was a latter-day astronomer-mathematician belonging to the Kerala school of astronomy and mathematics, gives the value of $\pi$ correct to 17 decimal places and the author of the book has claimed this value can be calculated using the series in Eq.\eqref{SankaraVarman} (see \cite{SankaraVarman} p.45). Also computations can be easily done by hand and only 20 terms are required to get an accuracy of 11 decimal places which is the accuracy attributed to M\=adhava.

\end{document}